\documentclass[11pt,a4paper]{article}
\usepackage[english]{babel}
\usepackage{amsmath, amsthm, amsfonts, amssymb, amscd}
\usepackage[utf8]{inputenc}
\usepackage[pdftex]{graphicx}
\usepackage[matrix,arrow]{xy}
\usepackage{enumerate}
\DeclareGraphicsExtensions{.bmp,.png,.pdf,.jpg}
\newtheorem*{proper}{Property}
\newtheorem{defi}{Definition}
\newtheorem{example}{Example}
\newtheorem{teo}{Theorem}
\newtheorem{lema}[teo]{Lemma}
\newtheorem{cor}[teo]{Corollary}
\newtheorem{prop}[teo]{Proposition}

\newtheorem{rem}{Remark}
\newcommand{\hookuparrow}{\mathrel{\rotatebox[origin=c]{90}{$\hookrightarrow$}}}

\DeclareMathOperator{\spec}{Spec}

\newcommand{\Oh}{\mathcal{O}}
\newcommand{\pf}{\noindent{\bf Proof\ \ }}
\newcommand{\cqd}{{\hfill $\rule{2mm}{2mm}$}\vspace{3mm}}
\newcommand{\Z}{\mathbb Z}
\newcommand{\N}{\mathbb N}
\newcommand{\tildeOh}{\widetilde{\mathcal{O}}}
\newcommand{\I}{\mathcal{I}}
\newcommand{\J}{\mathcal{J}}

\newcommand{\C}{\mathcal{C}}
\newcommand{\K}{\mathcal{K}}
\newcommand{\f}{\mathfrak f}
\usepackage{float}

\begin{document}
\title{Sets of Values of Fractional Ideals of Rings of Algebroid Curves}

\author{E. M. N. de Guzm\'an\footnote{Supported by a fellowship from CAPES} \\ {ICMC/USP} \and A. Hefez\footnote{Partially supported by the CNPq Grant 307873/2016-1} \\ {UFF}}

\date{ }
\maketitle

\noindent{\bf Abstract}
The aim of this work is to study sets of values of fractional ideals of rings of algebroid curves and explore more deeply the symmetry that exists among sets of values of dual pairs of  ideals when the ring is Gorenstein. We also express the codimension of a fractional ideal in terms of the maximal points of the value set of the ideal. Finally we also show that the Gorensteiness of a ring of an algebroid curve is equivalent to some conditions relating certain codimensions of fractional ideals and of their duals. \bigskip

\medskip

\section{Introduction}

Semigroup of values of rings of algebroid plane branches were studied by Zariski in \cite{Za1} and they constitute over $\mathbb C$ a complete set of discrete invariants for their topological classification. These semigroups live in $\N$ and they are determined by a finite set of elements, since, as it is very well known, they have a conductor, that is, a number $c$ such that $c-1$ is not in the semigroup, but from $c$ on all natural numbers belong to the semigroup. From the work of Ap\'ery \cite{Ap}, it follows also that this semigroup is symmetric, implying that in the interval $[0,c-1]$ there as many elements in the semigroup as in its complement. Many years later, E. Kunz, in \cite{Ku}, showed that the symmetry of the semigroup of values for branches in a higher dimensional space is equivalent to the property that the ring of the branch is Gorenstein.

For a plane curve singularity with several branches over $\mathbb C$, O. Zariski in \cite{Za1} characterized its topological type through the semigroups of its branches and their mutual intersection multiplicities. This characterization was shown by R. Waldi in \cite{W} to be equivalent to the knowledge of the semigroup of values of the curve, this time in $\N^r$, where $r$ is the number of branches of the curve. Although not finitely generated, this semigroup, for $r=2$, was shown by A. Garcia in \cite{Ga} to be determined in a combinatorial way by a finite set of points that he called {\em maximal points}. Garcia also showed that these maximal points of the semigroup of a plane curve have a certain symmetry. These results were generalized later for any value of $r$ by F. Delgado in \cite{D87}. Delgado discovered several types of maximals among them the relative and absolute maximals and showed that the relative maximals determine the semigroup of values in a combinatorial way and extended Garcia's symmetry showing that relative and absolute maximals determine each other. Some time later, Delgado, in \cite{D88}, generalizing the work of Kunz, introduced a concept of symmetry for the semigroup of values in $\N^r$ and showed that this symmetry is equivalent to the property of the ring of the curve to be Gorenstein.

Recently, D. Pol, in the work \cite{Po}, showed that the property of the singularity to be  Gorenstein is equivalent to the fact that the set of values of the dual of any fractional ideal is determined by the set of values of the ideal. This motivated us to describe the property of a curve singularity to be Gorenstein by means of Ap\'ery sets following the original steps by Ap\'ery. For this, we defined the Ap\'ery set for a set of values of a fractional ideal and characterized Gorenstein rings of algebroid curves by means of some symmetry property between the Ap\'ery set of a fractional ideal and that of its dual. 

On the other hand, the computation of the codimension of an ideal, using elementary paths in its set of values was also considered by Garcia and Delgado, being generalized for fractional ideals by D'Anna \cite{Da} and by Barucci, D'Anna and Fr\"oberg and in \cite{BDF1}. These authors observed that the basic properties of semigroup of values extend trivially to the sets of values of fractional ideals, implying that their relative maximal points  determine in the same combinatorial way such sets, as for semigroups of values. What were missing in the literature and we provide here is a formula relating the codimension of an ideal with the maximal points of its set of values. 

Finally we also show that the property of a ring of an algebroid curve to be Gorenstein is equivalent to some conditions relating codimensions of certain fractional ideals and of their duals. This will imply, for Gorenstein rings, a symmetry property between the absolute maximals of the set of values of a fractional ideal and the relative maximals of the set of values of its dual ideal, and vice versa, generalizing the result obtained by A. Campillo, F. Delgado and K. Kiyek in \cite{CDK} for semigroups of values. This last result was also recently obtained by D. Pohl in \cite{Po2} by other methods. 

\section{Algebroid curves}

This introductory section serves to set up our framework and it is based on the work of Gorenstein \cite{Go}, Hironaka \cite{Hi} and Garcia-Lax \cite{Ga-La}.

Throughout this work we will denote by $\Z$ the set of integers and by $\N$ the set of non-negative integers.

A reduced algebroid curve is given by $C=\spec(\Oh)$ for a finitely generated, local, reduced and complete $k$-algebra of pure dimension one $\Oh$, where $k$ is a coefficient algebraically closed field for $\Oh$. The ring $\Oh$ will be called the ring of $C$.

Let $\wp_1, \ldots,\wp_r$ be the minimal primes of $\Oh$. We will use the notation $I=\{1,\ldots,r\}$. The {\em branches} of $C$ are the integral schemes $C_i=\spec(\Oh_i)$, $i\in I$, where $\Oh_i=\Oh/\wp_i$. We will denote by $\pi_i\colon \Oh \to \Oh/\wp_i$ the canonical surjection, which corresponds to an inclusion $C_i \hookrightarrow C$. 

Since $\Oh$ is reduced, we have that $\bigcap_{i=1}^r \wp_i =\sqrt{(0)}=(0)$, so the homomorphism  
$$
\begin{array}{rcl}
\pi\colon \Oh & \hookrightarrow &\Oh_1\times \cdots \times \Oh_r\\
            h & \mapsto & (\pi_1(h),\ldots,\pi_r(h)) 
						\end{array}$$
is an injection.

More generally, if $J=\{i_1,\ldots,i_s\}$ is any subset of $I$, we may consider $C_J=\spec(\Oh_J)$, where 
$\Oh_J=\Oh/\cap_{j=1}^s \wp_{i_j}$. The scheme $C_J$ represents the curve $C_{i_1}\cup\cdots\cup C_{i_s}$. We will consider the natural surjection  
$$
\pi_J\colon \Oh \longrightarrow \Oh_J.
$$

Notice that we have $\Oh_I=\Oh$ and $C_I=C$.
\medskip

We will denote by $\mathcal K$ the total ring of fractions of $\Oh$ and by $\mathcal K_J$ the total ring of fractions of the ring
$\Oh_J$, when $J\subset I$. 

If $J=\{i\}$, then $\Oh_{\{i\}}$ is equal to the above defined domain $\Oh_i$ whose field of fractions will be denoted by $\mathcal K_i$. Let $\widetilde{\Oh}$ be the integral closure of $\Oh$ in $\mathcal K$ and $\widetilde{\Oh}_J$ be that of $\Oh_J$ in $\mathcal K_J$. By the splitting principle, one has that $\widetilde{\Oh}_J\simeq \widetilde{\Oh}_{i_1} \times \cdots \times \widetilde{\Oh}_{i_n}$, which in turn is the integral closure of $ {\Oh}_{i_1} \times \cdots \times {\Oh}_{i_n}$ in its total ring of fractions. \medskip

We have the following diagram:
\[
\begin{array}{rcc}
\mathcal K_J & \simeq & \mathcal K_{i_1} \times \cdots \times \mathcal K_{i_n} \\
\hookuparrow     &  & \hookuparrow \\
\widetilde{\Oh}_J &  \simeq & \widetilde{\Oh}_{i_1} \times \cdots \times \widetilde{\Oh}_{i_n}\\
\hookuparrow     &   & \hookuparrow\\
\Oh_J         & \hookrightarrow & \Oh_{i_1} \times \cdots \times \Oh_{i_n}	
	\end{array}
	\]
		
Since each $\widetilde{\Oh_i}$ is a complete DVR over the field $k$ with valuation $v_i$, one has that $\widetilde{\Oh}_i \simeq k[[t_i]]$, where $v_i(t_i)=1$. It follows that $\mathcal K_i\simeq k((t_i))$ is a valuated field with the extension of the valuation $v_i$. We also have that $\mathcal K \simeq k((t_1))\times \cdots \times k((t_r))$. \medskip

We will denote by $\varphi_i$ the map $\Oh \to \tildeOh_i$ and by $\varphi$ the injective map $\Oh \to  \tildeOh_1\times \cdots \times\tildeOh_r$, $\varphi(h)=(\varphi_1(h),\ldots,\varphi_r(h)).$\medskip

It is known that $\widetilde{\Oh}_J$ is a finitely generated $\Oh_J$-module, hence the set $$\mathcal C_J=\{h\in \mathcal K_J; \, h \widetilde{\Oh}_J \subset \Oh_J\}$$ has a regular element (ie. a nonzero divisor). It is easily verified that this set $\mathcal C_J$ is an ideal of both $\Oh_J$ and $\widetilde{\Oh}_J$ and, in fact, it contains any simultaneous ideal of these two rings. This ideal is called the {\em conductor ideal} of $\widetilde{\Oh}_J$ in $\Oh_J$. When $J=I$, the conductor ideal $\C_I$ will be denoted by $\C$.
\medskip

An important class of rings $\Oh$ are the {\em Gorenstein rings}, for which one has, by definition that, as $k$-vector spaces,
$$
\dim_k \Oh/\mathcal C= \dim_k \tildeOh/\Oh.
$$

The next result will characterize the ideals of $\widetilde{\Oh}$.

\begin{prop} If $\mathcal I$ is an ideal of $\widetilde{\Oh}\simeq k[[t_1]]\times \cdots \times k[[t_r]]$, then $\mathcal I= (t_1^{\gamma_1},\ldots,t_r^{\gamma_r})\widetilde{\Oh}$, where $\gamma_i =\min\{v_i(h); \ h\in \mathcal I\}$.
\end{prop}
\pf  We will prove that $(t_1^{\gamma_1},\ldots,t_r^{\gamma_r})$ generates $\mathcal I$. First of all let us observe that $(t_1^{\gamma_1},\ldots,t_r^{\gamma_r})\in \mathcal I$. Indeed, let $h^i\in \mathcal I$ be such that $v_i(\pi_i(h^i))=\gamma_i$, then if we take $\lambda_1,\ldots,\lambda_r \in k$ general enough, then $h=\lambda_1h^1+\cdots+\lambda_rh^r \in \mathcal I$ and $h=(t_1^{\gamma_1}u_1, \ldots, t_r^{\gamma_r}u_r)$, where each $u_i$ is a unit  in $k[[t_i]]$ , hence $h=u(t_1^{\gamma_1}, \ldots, t_r^{\gamma_r})$, where $u$ is a unit in $\widetilde{\Oh}$. This shows that $$(t_1^{\gamma_1}, \ldots, t_r^{\gamma_r})=u^{-1}h \in \mathcal I.$$

Now, any element $g\in \mathcal I$ may be written as
\[
g=(t_1^{\gamma_1+\alpha_1}u_1, \ldots, t_r^{\gamma_r+\alpha_r}u_r)=(t_1^{\alpha_1}u_1, \ldots, t_r^{\alpha_r}u_r)(t_1^{\gamma_1}, \ldots, t_r^{\gamma_r}).
\]
\cqd

It then follows that there is some $\gamma =(\gamma_1,\ldots,\gamma_r)\in \N^r$ such that $\mathcal C= \underline{t}^\gamma\widetilde{\Oh}$, where  $\underline{t}^\gamma=(t_1^{\gamma_1}, \ldots, t_r^{\gamma_r})$.\medskip
	
\begin{example}{\rm(Ap\'ery \cite{Ap})} If	$\spec(\Oh)$ is a plane branch, then one has that $\mathcal C=t^{2\delta}k[[t]]$, where $\delta=\dim_k\frac{\widetilde{\Oh}}{\Oh}=\dim_k\frac{\Oh}{\mathcal C}$. The number $c=2\delta$ is the degree of the conductor $\mathcal C$.\end{example}

\begin{example}{\rm (Gorenstein \cite{Go})} If $\spec(\Oh)$ is a reduced plane curve given by an equation $f=f_1\cdots f_r$, then the conductor of $\widetilde{\Oh}$ in $\Oh$ is a principal ideal generated by $(t_1^{I_1+c_1}, \ldots, t_r^{I_r+c_r})$, where the $c_i$ denote the degrees of the conductors ideals $\C_i$ of the branches determined by the $f_i$ and  $I_i=\sum_{j\neq i}I(f_j,f_i)$. \end{example}

These two examples will be stated in greater generality later in Theorems 6 and 7.

Let us denote by $Z(A)$ the set of zero divisor of a ring $A$. One has that $Z(\Oh)=\cup_i\wp_i$ and if $A=\widetilde{\Oh}$ or $A=\mathcal K$, then $Z(A)=\{(a_1,\ldots,a_r)\in A; \, a_i=0, \ \text{for \ some} \ i\}$. So, one has that $g\in Z(\Oh)$ (or equivalently $\varphi_i(g)=0$), if and only if $\pi_i(g)=0$, for some $i=1,\ldots,r$.\bigskip

\begin{defi} \label{intersection} If $J_1$ and $J_2$ are two disjoint subset of $I$, then the intersection multiplicity of $C_{J_1}$ and $C_{J_2}$ is defined by
$$
I(C_{J_1},C_{J_2})=\dim_k \frac{\Oh}{\bigcap_{j\in J_1}\wp_j+\bigcap_{j\in J_2}\wp_j}.
$$
\end{defi}

In the above situation we also have that $C_{J_1}\cup C_{J_2}=C_{J_1\cup J_2}$.
    
\section{Fractional ideals}

An $\Oh$-submodule $\mathcal I$ of $\mathcal K$ will be called a {\em fractional ideal} of $\Oh$ if it  contains a regular element of $\Oh$ and there is a regular element $d$ in $\Oh$ such that $d\,\mathcal I \subset \Oh$. Such an element $d$ will be called a common denominator for $\mathcal I$.

Since $d\,\mathcal I$ is an ideal of $\Oh$, which is a noetherian ring, one has that $\mathcal I\subset \mathcal K$ is a nontrivial fractional ideal if and only if it contains a regular element of $\Oh$ and it is a finitely generated $\Oh$-module.

Examples of fractional ideals of $\Oh$ are $\Oh$ itself, $\widetilde{\Oh}$, the conductor $\mathcal C$ of $\tildeOh$ in $\Oh$, or any ideal of $\Oh$ or of $\tildeOh$ that contains a regular element. Also, if $\I$ is a fractional ideal of $\Oh$, then for all $\emptyset \neq J\subset I$ one has that $\pi_J(\I)$ is a fractional ideal of $\Oh_J$.

\begin{rem}
Any element $h\in \mathcal K$ will be written as $h=\frac{h_1}{h_2}$, where $h_1,h_2\in \Oh$, $h_2\notin Z(\Oh)$. If $h=\frac{h_1}{h_2}\in \mathcal I$, then $h_1\in \mathcal I$. And, $h\in Z(\mathcal K)$ if and only if $h_1\in Z(\Oh)$. 
\end{rem}

So, we may assume that the regular element contained in a fractional ideal is indeed in $\Oh$. In the set of fractional ideals of $\Oh$ there are defined the operations $+$, $\cdot$ and $\colon$. The first two operations are clear and the last one is defined, for $\mathcal I \neq \{0\}$, as usual:
\[
(\mathcal I\colon \mathcal J)=\{x\in \mathcal K; \, x\mathcal  J \subset \mathcal I\}.
\]

This defines a fractional ideal. Indeed, it is an $\Oh$-module which contains a regular element: $dd'$, where $d$ is a common denominator for $\mathcal I$ and $d'$ is a common denominator for $\mathcal J$. 
On the other hand, the element $b d$, where $b$ is a regular element in $\mathcal J \cap \Oh$, is a common denominator for $(\mathcal I\colon \mathcal J)$.\medskip

As an example of this construction we have that $\mathcal C=(\Oh:\tildeOh)$.

\begin{prop} Given a fractional ideal $\mathcal I$ of $\Oh$ there exist $\beta\in \N^r$ and $\alpha \in \Z^r$ such that
\[
\underline{t}^\beta \tildeOh \subset \mathcal I \subset  \underline{t}^\alpha \tildeOh.
\]
\end{prop}
\pf \noindent From the definition of a fractional ideal we know that there exists a regular element $d$ in $\Oh$ such that $d\,\mathcal I \subset \Oh \subset \tildeOh$. Therefore, as subsets of $\mathcal K$, one has that $\mathcal I \subset \frac{1}{d} \tildeOh$. Since we may write
$\frac{1}{d}=\underline{t}^\alpha u$, where $\alpha\in \Z^r$ and $u$ is a unit of $\tildeOh$, we get that
\[
\mathcal I \subset \underline{t}^\alpha u\tildeOh=\underline{t}^\alpha \tildeOh.
\]

On the other hand, we now that there is a regular element of $\Oh$ in $\mathcal I$. Write $a=\underline{t}^\epsilon v$, where $\epsilon \in \N^r$ and $v$ is a unit of $\tildeOh$. Let $\underline{t}^\gamma$ be a generator of the conductor ideal of $\tildeOh$ in $\Oh$, then 
\[
\underline{t}^{\epsilon+\gamma} \tildeOh =  \underline{t}^\epsilon  v \underline{t}^\gamma   \tildeOh \subset a\Oh\subset \mathcal I.
\]
The remaining inclusion we are looking for follows taking $\beta=\epsilon +\gamma$.
\cqd

We define the {\em dual module} of $\mathcal I$ as being the $\Oh$-module $\mathcal I^{\vee} =Hom_\Oh(\mathcal I, \Oh)$.\medskip

If $\varphi\in \mathcal I^{\vee}$ and $h=\frac{h_1}{h_2}\in \mathcal I$, then it is easy to verify that 
$$ \varphi\left(\frac{h_1}{h_2}\right)=\frac{\varphi(h_1)}{h_2}.
$$

Our aim now is to realize  $\mathcal I^{\vee}$ as a fractional ideal of $\Oh$. To do this, we begin with the following result.

\begin{lema} Let $\varphi \in \mathcal I^{\vee}$.
If $g \in \mathcal I \setminus Z(\mathcal K)$, then $\frac{\varphi(g)}{g}$ is independent of $g$.
\end{lema}
\noindent \pf  Let $g=\frac{g_1}{g_2}$ and $g'=\frac{g_1'}{g_2'}$ in $\mathcal I \setminus Z(\mathcal K)$, then
$$\begin{array}{l}
\frac{\varphi(g)}{g}=\frac{g_2\varphi(g)}{g_1}=\frac{\varphi(g_2g)}{g_1}=\frac{\varphi(g_1)}{g_1},\\ \\
\frac{\varphi(g')}{g'}=\frac{g_2'\varphi(g')}{g_1'}=\frac{\varphi(g_2'g')}{g_1'}=\frac{\varphi(g_1')}{g_1'},
\end{array}
$$
which make sense since $g_1$ and $g_1'$ are not in $Z(\Oh)$.

But, since $\varphi(g_1g_1')=\varphi(g_1'g_1)$ and $\varphi$ is an $\Oh$-modules homomorphism, it follows that $g_1\varphi(g_1')=g_1'\varphi(g_1)$, which implies our result.
\cqd

For $g\in \mathcal I \setminus Z(\mathcal K)$, the map
$$
\begin{array}{rcl}
\Phi\colon Hom_\Oh(\mathcal I, \Oh) & \to & \mathcal K\\
\varphi & \mapsto & \frac{\varphi(g)}{g},
\end{array}
$$
which is independent of $g$, is a homomorphism of $\Oh$-modules.\medskip

Given a fractional ideal $\mathcal I$ of $\Oh$, we will use the notation
\[
\mathcal I^*=(\Oh\colon \mathcal I)=\{x\in \mathcal K; \, x\mathcal I \subset \Oh\}.
\]

\begin{prop} The map $\Phi$ induces an $\Oh$-module isomorphism between $\mathcal I^{\vee}$ and $\mathcal I^*$.
\end{prop}
\noindent \pf We first show that $\Phi$ is injective. Let $\varphi\in Hom_\Oh(\mathcal I, \Oh)$ be such that $\frac{\varphi(g)}{g}=0$, for some $g\in (\mathcal I\cap \Oh) \setminus Z(\Oh)$. Let $h_1\in \mathcal I \cap \Oh$, then
$$
\varphi(h_1)=\frac{g\varphi(h_1)}{g}=h_1 \frac{\varphi(g)}{g}=0,
$$
which implies that $\varphi(\frac{h_1}{h_2})=\frac{\varphi(h_1)}{h_2}=0$, for $\frac{h_1}{h_2}\in \mathcal I$, so $\Phi$ is injective.

Now we show that $Im(\Phi) \subset  \mathcal I^*$. Indeed, let $a\in Im(\Phi)$, then $a=\frac{\varphi(g)}{g}$ for some $\varphi \in \mathcal I^*$ and some $g\in (\mathcal I \cap \Oh) \setminus Z(\Oh)$. Now, let $h=\frac{h_1}{h_2} \in \mathcal I$, then 
$$
a\frac{h_1}{h_2}=\frac{\varphi(g)}{g} \frac{h_1}{h_2} =\frac{\varphi(gh_1)}{gh_2}=\frac{g}{gh_2}\varphi(h_1) =\frac{\varphi(h_1)}{h_2}= \varphi\left(\frac{h_1}{h_2}\right)\in \Oh.
$$

Finally, let us show that $\mathcal I^* \subset Im(\Phi)$. Take $a$ in $\mathcal I^* $ and define 
$$
\begin{array}{rcl}
\varphi_a\colon \mathcal I & \to & \Oh\\
h& \mapsto & ah.
\end{array}
$$
It is easy to verify that $\varphi_a$ is a homomorphism of $\Oh$-modules, and 
$$a=\frac{ag}{g}=\frac{\varphi_a(g)}{g} \in Im(\Phi).
$$
\cqd

We will always view $\mathcal I^{\vee}=Hom_\Oh(\mathcal I, \Oh)$, as the fractional ideal $\mathcal I^*$, which will be called the dual ideal of the fractional ideal $\mathcal I$. It is easy to verify that taking duals is order reversing with respect to inclusions; that is,
\[
\mathcal I \subset \mathcal J \, \Longrightarrow \, \mathcal I^* \supset \mathcal J^*.
\]

As examples of duals of fractional ideals one has that $\Oh^*=\Oh$ and $\tildeOh^*=\mathcal C$. 
\medskip

One always has an inclusion $\mathcal I \subset ({\mathcal I}^*)^*$, which is almost tautological. On the other hand, it is not always true that $({\mathcal I}^*)^* \subset \mathcal I$; but one has the following result. 

\begin{teo}{\rm (\cite[Theorem 6.3]{B})} The following statemens are equivalent:
\begin{enumerate}[i)]
\item $\Oh$ is a Gorenstein ring;
\item $({\mathcal I}^*)^* \subset \mathcal I$ for all fractional ideal $\mathcal I$ of $\Oh$;
\item $\dim_k \frac{\mathcal J}{\mathcal I}=\dim_k \frac{\mathcal I^*}{\mathcal J^*}$, for all fractional ideals $\mathcal I$ and 
 $\mathcal J$ such that $\mathcal I \subset \mathcal J$.  
\end{enumerate}
\end{teo}\smallskip

If $J\subset I$ we define
$$
\delta_J= \dim_k \widetilde{\Oh_J}/\Oh_J,
$$
and if we put $J'=I\setminus J$, then one defines $I_J= I(C_J,C_{J'})$, we have the following results:

\begin{teo}{\rm (\cite[Corollary 3.10]{Ga-La})} Let $\Oh$ be a Gorenstein ring, $J\subset I$  and $T_1,\ldots , T_t$ a partition of $I$. Then one has
$$ 
\delta_J+\frac 12 I_J=\sum_{j=1}^t\delta_{T_j}+\frac 12 \sum_{j=1}^tI_{T_j}.
$$
\end{teo}

\begin{teo}{\rm (\cite[Theorem 3.13]{Ga-La})}\label{conductor} Let $\Oh$ be a Gorenstein ring, then the conductor $\mathcal C$ of $\tildeOh$ in $\Oh$ is given by
$$
\mathcal C= (t_1^{I_1+2\delta_1}, \ldots, t_r^{I_r+2\delta_r})\tildeOh.
$$
\end{teo}

\section{Value sets of fractional ideals}

For any $n$, we will consider on $\Z^n$ the partial order given by
\[
(a_1,\ldots,a_n)\leq (b_1,\ldots,b_n) \ \Longleftrightarrow \ a_i \leq b_i, \ i=1,\ldots,n.
\]

\noindent We will write $(a_1,\ldots,a_n) < (b_1,\ldots,b_n)$ when $a_i<b_i$, for all $ i=1,\ldots,n$.\medskip

Consider the value map
\[
\begin{array}{ccc}
v \colon \mathcal K \setminus Z(\mathcal K) & \to &  \Z^r\\
   h & \mapsto &  (v_1(\pi_1(h)), \ldots ,v_r(\pi_r(h))),
	\end{array}
\]
where $\pi_i$ here denotes the projection $\mathcal K \to \mathcal K_i$, which is the extension of the previously defined projection map $\pi_i\colon \Oh \to \Oh_i$. \medskip
		
If $\mathcal I$ is a fractional ideal of $\Oh$, we define the {\em value set} of $\mathcal I$ as being 
$$E=v(\mathcal I\setminus Z(\mathcal K))\subset \Z^r.$$

We will denote by $E^*$ the value set of $\mathcal I^*$. The value set of $\Oh$ will be denoted by $S(\Oh)$, or simply by $S$,  and it is a subsemigroup of $\N^r$, called the semigroup of values of $\Oh$. The value sets $E$ of fractional ideals are not necessarily closed under addition, but they are such that $S+E\subset E$. For this reason they are called monomodules over the semigroup $S$.

From the definition of the ideal $\I^*$, if $E$ is the value set of $\I$ and $E^*$ is the value set of $\I^*$, it follows immediately that 
\begin{equation}\label{$E+E^*$}
E+E^*\subset S(\Oh).
\end{equation}

Let $I=\{1,\ldots, r\}$ and $J\subset I$. If $J=\{i_1,\ldots,i_n\}$, then we denote by
$\pi_J$ the projection $\mathcal K \to \mathcal K_{i_1} \times \cdots \times \mathcal K_{i_n}$ and by $pr_J$ the corresponding projection $\Z^r \to \Z^n$, $(\alpha_1,\ldots,\alpha_r) \mapsto (\alpha_{i_1},\ldots,\alpha_{i_n})$. \medskip

Let us define 
$$E_J=v(\pi_J(\I)\setminus Z(\K_J)). 
$$

\noindent {\bf Caution} \ If $j\in J=\{i_1,\ldots,i_t,\ldots i_s\} \subset I$, with $i_t=j$, for $\beta=(\beta_{i_1},\ldots,\beta_{i_s})\in E_J$, we will define $$pr_j(\beta)=\beta_{i_t}=\beta_j,$$ instead of $pr_j(\beta)=\beta_{i_j}$, as it would be natural.\medskip

From Proposition 2 one has that, given a set of value $E$, there exist $\beta\in E$ and $\alpha \in \Z^r$ such that 
\[\beta +\N^r \subset E \subset \alpha +\N^r.
\]

Value sets of fractional ideals have the following fundamental properties, analogous to the properties of semigroups of values described by Garcia for $r=2$ in \cite{Ga} and by Delgado for $r\geq 2$ in \cite{D87} (see also \cite{Da} or \cite{BDF1}):

\begin{proper}[\textbf{A}]
	If $\alpha=(\alpha_1,\ldots,\alpha_r)$ and $\beta=(\beta_1,\ldots,\beta_r)$ belong to $E$, then 
	$$min(\alpha,\beta)=(min(\alpha_1,\beta_1),\ldots,min(\alpha_r,\beta_r))\in E.$$
\end{proper}

\begin{proper}[\textbf{B}]
	If $\alpha=(\alpha_1,\ldots,\alpha_r), \beta=(\beta_1,\ldots,\beta_r)$ belong to $E$, $\alpha\neq\beta$ and $\alpha_i=\beta_i$ for some $i\in\{1,\ldots,r\}$, then there exists $\gamma\in E$ such that $\gamma_i>\alpha_i=\beta_i$ and $\gamma_j\geq min\{\alpha_j,\beta_j\}$ for each $j\neq i$, with equality holding if $\alpha_j\neq\beta_j$.
\end{proper}

Since there is an $\alpha \in \Z^r$ such that $E \subset \alpha +\N^r$, then from Property (A) there exists a unique $m=m_E=\min(E)$, that is, for all $(\beta_1,\ldots,\beta_r)\in E$, one has $\beta_i\geq m_i$, $i=1,\ldots,r$.\medskip

On the other side, one has the following
\begin{lema} If $\beta, \beta'\in E$ are such that  $\beta +\N^r \subset E$ and $\beta' +\N^r \subset E$, then
$\min(\beta,\beta') +\N^r \subset E$.
\end{lema}
\pf
Let $a\in \N^r$ be arbitrary. Since  $\beta +a\in E$ and $\beta'+a\in E$, from Property (A), we have that 
$$
\min(\beta,\beta')+a =\min(\beta+a,\beta'+a)\in E.
$$

Therefore, $min(\beta,\beta')+\mathbb{N}^r\subset E$.
\cqd

The above lemma guarantees that there is a unique least element $\gamma\in E$ with the property that $\gamma+\N^r \subset E$. This element is what we call the conductor of $E$ and will denote it by $c(E)$.

Observe that one always have 
\[
c(E_J)\leq c(pr_J(E)), \ \ \forall\, J\subset I.
\]

One has the following result:

\begin{lema} \label{proj} If $\I$ be a fractional ideal of $\Oh$ and $J\subset I$, then
$pr_J(E)=E_J$.
\end{lema}
\pf One has obviously that $pr_J(E)\subset E_J$. On the other hand, let $\alpha_J \in E_J$. Take $h\in \I$ such that $v_J(\pi_J(h))=\alpha_J$. If $h\not\in Z(\K)$ we are done. Otherwise, choose any $h'\in \I \setminus Z(\K)$ such that $pr_J(v(h'))>\alpha_J$, which exists since $E$ has a conductor. Hence, $v_J(h+h')=\alpha_J$, proving the other inclusion.
\cqd

Let $S=S(\Oh)$. We will use the following notation:
\[\mathfrak f(S) =c(S)-(1,\ldots,1),\]
and call $\f(S)$ the {\em Frobenius vector} of $S$ or of $\Oh$.\medskip

\section{Maximal points} \label{maximals}

We now introduce the important notion of a {\em fiber} of an element $\alpha\in E$ with respect to a subset $J\subset I=\{1,\ldots,r\}$ that will play a central role in what follows.

\begin{defi} 	Given $\alpha\in\mathbb{Z}^r$ and $\emptyset \neq J\subset I$, we define: 	
$$F_J(\alpha)=\{\beta\in\mathbb{Z}^r; pr_J(\beta)=pr_J(\alpha) \ \text{and} \ pr_{I\setminus J}(\beta)>pr_{I\setminus J}(\alpha)\},$$
$$\overline{F}_J(\alpha)=\{\beta\in\mathbb{Z}^r; pr_J(\beta)=pr_J(\alpha), \  \text{and} \ pr_{I\setminus J}(\beta) \geq pr_{I\setminus J}(\alpha)\},$$
$$F_J(E,\alpha)=F_J(\alpha)\cap E  \quad \text{and} \quad \overline{F}_J(E,\alpha)=\overline{F}_J(\alpha)\cap E,$$
$$F(\alpha)=\bigcup_{i=1}^rF_{\{i\}}(\alpha), \quad F(E,\alpha)=F(\alpha)\cap E.$$
	
The last set above will be called the {\em fiber} of $\alpha$.
	
\end{defi}

The sets $F_{\{i\}}(\alpha)$ and $\overline{F}_{\{i\}}(\alpha)$ will be denoted simply by $F_{i}(\alpha)$ and $\overline{F}_{i}(\alpha)$. Notice that $F_I(\alpha)=\{\alpha\}$.

\begin{defi}
Let $\alpha\in E$. We will say that $\alpha$ is a\textbf{ maximal} point of $E$ if $F(E,\alpha)=\emptyset$. 
\end{defi}

This means that there is no element in $E$ with one coordinate equal to the corresponding coordinate of $\alpha$ and the other ones bigger.	\medskip

By convention, when $r=1$, the maximal points of $E$ are exactly the gaps of $E$.

From the fact that $E$ has a minimum $m_E$ and a conductor $\gamma=c(E)$, one has immediately that all maximal elements of $E$ are in the limited region
\[
\{(x_1,\ldots,x_r)\in \Z^r; \ m_{Ei}\leq x_i < \gamma_i, \ \ i=1,\ldots,r\}.
\]

This implies that $E$ has finitely many maximal points.

\begin{defi} We will say that a maximal point $\alpha$ of $E$ is an \textbf{absolute maximal} if $F_J(E,\alpha)=\emptyset$ for every $J\subset I$, $J\neq I$. If a maximal point $\alpha$ of $E$ is such that $F_J(E,\alpha)\neq\emptyset$, for every $J\subset I$ with $\#J\geq2$, then $\alpha$ will be said to be a \textbf{relative maximal} of $E$.
\end{defi}

\begin{figure}[h]
 \centering
\includegraphics[scale=1.1]{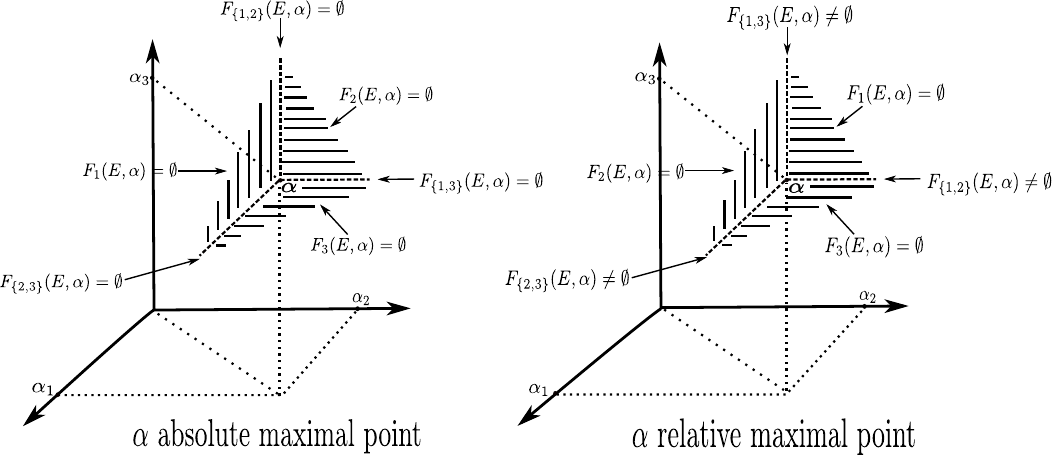}
\caption{Maximal points}
\end{figure}

Notice that when $r=1$, as we said before, the maximals coincide with the gaps

In the case where $r=2$ the notions of maximal, relative maximal and absolute maximal coincide. For $r=3$ we only may have relative maximals or absolute maximals.\medskip

We will denote by $M(E)$, $RM(E)$ and $AM(E)$ the sets of maximals, of relative maximals and absolute maximals of the set $E$, respectively.\medskip

The theorem below says that the set of relative maximal of $E$ determines $E$ in a combinatorial sense as follows:
\begin{teo}[generation]\label{generation}
	Let $\alpha\in\mathbb{Z}^r$ be such that $p_J(\alpha)\in E_J$ for all $J\subset I$ with $\#J=r-1$. Then
	$$\alpha\in E \Longleftrightarrow \alpha\notin F(\mathbb{Z}^r,\beta), \ \forall \beta\in RM(E).$$
\end{teo}

We will omit the proof since this result is a slight modification of \cite[Theorem 1.5 ]{D87} with essentially the same proof.\medskip

The following result is the content of \cite[Corollary 1.9 and Corollary 2.7]{D88}.

\begin{teo} \label{D2.7} Let $\Oh$ be the ring of an algebroid curve and let $S=S(\Oh)$. One has
\begin{enumerate}[\rm i)]
\item $F(S,\f(S))=\emptyset$;
\item If $\Oh$ is Gorenstein, then $\f(S)$ is a relative maximal of $S$.
\end{enumerate}
\end{teo}

Part (i) of the above theorem implies the following result:

\begin{cor} \label{halfsymm} Let $\Oh$ be the ring of an algebroid curve, and let $S=S(\Oh)$. Then 
\[
\alpha \in S \Longrightarrow \ F(S, \f(S)-\alpha)=\emptyset.
\]
\end{cor}
\pf
Indeed, if $\alpha\in S$ and $ F(S, \f(S)-\alpha)\neq\emptyset$, then if $\theta\in  F(S, \f(S)-\alpha)$, then  $\theta+\alpha \in  F(S, \f(S))$, which is a contradiction. 
\cqd

\section{Symmetry}
In what follows we will denote by $\Oh$ the ring of an algebroid curve and by $S=S(\Oh)$ its semigroup of values. Recall that we defined $\f(S)$ as being $c(S)-(1,\ldots,1)$. A fractional ideal of $\Oh$ will be denoted by $\I$, its dual by $\I^*$, while  $E$ and  $E^*$ will denote the value sets of $\I$ and $\I^*$, respectively.\medskip

The property of the ring $\Oh$ of an algebroid curve to be Gorenstein is equivalent to the symmetry of the semigroup of values $S(\Oh)$ of $\Oh$ as pointed out by E. Kunz in \cite{Ku} in the irreducible case and generalized by F. Delgado de la Mata in \cite{D88} as follows:

\begin{teo}{\rm (\cite[Theorem 2.8]{D88})}\label{delgado88} The ring $\Oh$ is Gorenstein if and only if 
\begin{equation}\label{symmetry}
\forall \alpha\in \N^r, \ \ \alpha \in S \ \ \Longleftrightarrow \ \ F(S,\f(S) -\alpha)=\emptyset.
\end{equation}
\end{teo}

A semigroup $S$ satisfying the property (\ref{symmetry}) above, will be called a {\em symmetric semigroup}.\medskip

The above result was generalized by D. Pol in \cite{Po} as follows:

\begin{teo}{\rm (\cite[Theorem 2.4]{Po})} \label{poltheo} The ring $\mathcal{O}$ is Gorenstein if and only if,  for all fractional ideal $\mathcal{I}$ of $\mathcal{O}$, the following property is verified 
	\begin{equation}\label{pol}
	\forall \beta\in \Z^r, \ \ \beta\in E^* \ \ \Longleftrightarrow \ \ F(E,\f(S)-\beta)=\emptyset.
	\end{equation}
	\end{teo}
	
Property (\ref{pol}) says that for a Gorenstein ring $\Oh$, the values set $E$ of any fractional ideal,  determines the values set $E^*$ of its dual $\I^*$, and conversely.\medskip 

Let us recall that in the case of a plane branch, the Gorenstein property of its ring was a consequence of the symmetry of its semigroup of values. This symmetry was proved by Ap\'ery in \cite{Ap} by showing that in this case the nowadays called Ap\'ery sets have a very special property. We will generalize this result to any algebroid curve and for this, we need to generalize the concept of Ap\'ery sets.  

Let $\I$ be a fractional ideal of $\mathcal{O}$ and let $E$ be its set of values. For any $\alpha\in E\setminus \{0\}$, we define the \emph{Ap\'ery set} of $E$ \emph{with respect to} $\alpha$ to be the set
$$A_\alpha(E)=\{\beta\in E; \ \beta-\alpha\notin E\}.$$
We will denote $A_\alpha(S)$ by $A_\alpha$.

\begin{lema}\label{Ap} Let $S$ be the set of values of an algebroid curve, let $\alpha\in S\setminus \{0\}$ and let $A_\alpha$ be its corresponding Ap\'ery set. For $\beta\in \mathbb{Z}^r$ one has that  $\beta\notin S$ if only if $\beta=a-\rho\alpha$, for some $a\in A_\alpha$ and some $\rho\in\N\setminus \{0\}$.
\end{lema}
\pf

\noindent ($\Rightarrow)$ Since $S\subset \N^r$, we have that $\alpha >0$. Suppose that $\beta\notin S$ and let $\rho$ be the positive integer such that $a=\beta+\rho\alpha\in S$, but  $\beta+(\rho-1)\alpha\notin S$. This implies $a\in A_\alpha$, giving the result.

\noindent ($\Leftarrow)$ We have that $\beta=a-\rho\alpha$ for $a\in A_\alpha$ and $\rho>0$. Suppose that $\beta\in E$. Since $\alpha \in S$ and $\rho>0$, we have that 
$$a-\alpha=\beta+ (\rho-1)\alpha  \in S,$$
which contradicts the fact that $a\in A_\alpha$.
\cqd

Now, we can state the following result.

\begin{teo}\label{ed}
Let $S=S(\Oh)$ be a semigroup of values in $\N^r$ and let $\alpha\in S\setminus \{0\}$. The following conditions are equivalent:
\begin{enumerate}[\rm i)]
\item $\Oh$ is Gorenstein;
\item $S$ is symmetric;
\item $\forall a \in \N^r, \ a\in A_\alpha \ \ \Longrightarrow \ \ \emptyset \neq F(S,\f(S)+\alpha-a)\subseteq A_\alpha$.
\end{enumerate}
\end{teo}
\pf  \  ${\rm i) \Leftrightarrow ii)}$ This is the content of Theorem \ref{delgado88}.\medskip

\noindent ${\rm ii) \Rightarrow iii)}$ Suppose that $S$  is symmetric. Let $a\in A_\alpha$, then $a-\alpha\notin S$ and since $S$ is symmetric, we have $F(S,\f(S)-(a-\alpha))\neq\emptyset$. Now we will show that $F(S,\f(S)-(a-\alpha))\subseteq A_\alpha$.

Let $\gamma$ be in  $F(S,\f(S)-(a-\alpha))$. If $\gamma-\alpha$ belonged to $S$, then we would have $\gamma-\alpha\in F(S,\f(S)-a)$, then $F(S,\f(S)-a)\neq \emptyset$. But, since $a\in A_\alpha\subset S$ and $S$ is symmetric, it follows that $F(S,\f(S)-a)=\emptyset$, which is a contradiction. So, $\gamma-\alpha\notin S$ and consequently $\gamma \in  A_\alpha$. \medskip
	

\noindent ${\rm iii) \Rightarrow ii)}$	Assuming (iii), let us prove that $S$ is symmetric, i.e., $\gamma\in S \ \Leftrightarrow \ F(S,\f(S)-\gamma)=\emptyset$. 

We know from Corollary \ref{halfsymm} that the implication  $\gamma\in S \ \Rightarrow \ F(S,\f(S)-\gamma)=\emptyset$ is always satisfied. 

To prove the other direction, we show that if $\gamma\notin S$ then  $F(S,\f(S)-\gamma)\neq\emptyset$. So, let $\gamma\notin S$, then by Lemma \ref{Ap} we have that $\gamma=a+\rho\alpha$, for some $a\in A_\alpha$ and some negative integer $\rho$. From our hypothesis we know that $ F(S,\f(S)+\alpha-a)\neq \emptyset$, hence take $\theta$ in this set. 

Since $\theta \in S$ and $-(\rho+1)\alpha \in S$, we have that
$$\theta-(\rho+1)\alpha\in F(\mathbb{N}^\rho,\f(S)-(a+\rho\alpha))\cap S=F(S,\f(S)-\gamma),$$
showing that $F(S,\f(S)-\gamma)\neq\emptyset$ and therefore $S$ is symmetric.
\cqd

Observe that the condition $F(S,\f(S)+\alpha-a)\neq \emptyset$ implies that $a-\alpha \notin S$, hence if $a\in S$, then $a\in A_\alpha$. It follows that condition (iii) in the above proposition may be readed as 
\[
\forall a \in \N^r, \ a\in A_\alpha \ \ \Longleftrightarrow \ \ a\in S \ \text{and} \ \emptyset \neq F(S,\f(S)+\alpha-a)\subseteq A_\alpha.
\]
\begin{teo}\label{ed2}
	Let $\mathcal{O}$ be the ring of an algebroid curve. We denote by $\I$ a fractional ideal of $\Oh$ and by $E$ its value set. For $\alpha\in E\cap E^*$, the following conditions are equivalent:
	\begin{enumerate}[\rm i)]
	\item $\Oh$ is Gorenstein;
	\item  $\forall \, \I, \ \forall\, \beta\in \Z^r, \ \  \beta\in E^*\Longleftrightarrow F(E,\f(S)-\beta)=\emptyset$.
	\item $\forall \, \I,   \ \forall\,  a\in \Z^r, \ \ a\in A_\alpha(E) \ \Longrightarrow \ \emptyset \neq F(E^*,\f(S)+\alpha-a)\subseteq A_\alpha(E^*)$.
	\end{enumerate} 
\end{teo}
\pf ${\rm i) \Longleftrightarrow ii)}$ This is the content of Theorem \ref{poltheo}.\medskip

\noindent ${\rm ii) \Rightarrow iii)}$ \  Let us assume that (ii) is fulfilled. 
Let $a\in A_\alpha(E)$, then $a-\alpha\notin E$, so from (ii), we have $F(E^*,\f(S) +\alpha-a)\neq\emptyset$ (recall that (ii) is equivalent to $\Oh$ is Gorenstein and, in this situation, $(E^*)^*=E$). 

Now, suppose that  $\gamma\in F(E^*,\f(S)+ \alpha-a)$, then $\gamma \in E^*$. It remains to show that $\gamma-\alpha\notin E^*$. Suppose that the opposite holds, that is, $\gamma-\alpha\in E^*$, then $\gamma-\alpha\in F(E^*,\f(S)-a)$, hence $F(E^*,\f(S)-a)\neq \emptyset$. But, since $a\in E$, from (ii) one gets that $F(E^*,\f(S)-a)=\emptyset$, a contradiction.
\medskip

\noindent ${\rm iii)\Rightarrow i)}$ \ Let us assume that (iii) holds for all fractional ideal $\I$. In particular it holds for $\I=\Oh$. Now, Theorem \ref{ed} implies that $\Oh$ is Gorenstein.
\cqd

\section{Colengths of fractional ideals}

Let $\Oh$ be the ring of an algebroid curve and $\J\subset \I$ are two fractional ideals with sets of values $D$ and $E$, respectively. Since $\J\subset \I$, one has that $D\subset E$, hence $c(E)\leq c(D)$. Our aim in this section is to find a formula for the length $\ell_\Oh(\I/\J)$ of $\I/\J$ as $\Oh$-modules, called the colength of $\J$ in $\I$, in terms of the value sets $D$ and $E$. 

The motivation comes from the case in which $\Oh$ is a domain, when one has that
\[
\ell_\Oh(\I/\J)=\#(E\setminus D).
\]\smallskip

For $\alpha\in \Z^r$ and $\I$ a fractional ideal of $\Oh$, with value set $E$, we define 
\[
\I(\alpha)=\{h\in \I\setminus Z(\K); \ v(h)\geq \alpha\}.
\]

It is clear that if $m_E=\min E$, then $\I(m_E)=\I$.\medskip

One has the following result:

\begin{prop}{\rm (\cite[Proposition 2.7]{BDF1})}\label{prop2}
Let $\J\subseteq\I$ be two fractional ideals of $\Oh$, with value sets $D$ and $E$, respectively, then  
$$\ell_\Oh\left(\dfrac{\I}{\J}\right) = \ell_\Oh\left( \dfrac{\I}{\I(\gamma)}\right)-\ell_\Oh\left(\dfrac{\J}{\J(\gamma)}\right),$$ 
for sufficiently large $\gamma\in \N^r$ (for instance, if $\gamma \geq c(D)$).
\end{prop}

Let $e_i\in \Z^r$ denote the vector with zero entries except the $i$-th entry which is equal to $1$, and recall the definition we gave in Section 5 for $\overline{F}_i(E,\alpha)$,  then the following result will give us an effective way to calculate colengths of ideals.

\begin{prop}{\rm \cite[Proposition 2.2]{Da}} \label{prop1}
If $\alpha\in\mathbb{Z}^r$, then we have
$$\ell_\Oh \left(\dfrac{\I(\alpha)}{\I(\alpha+e_i)}\right)=\left\{
\begin{array}{ll}
1, & \ if \  \overline{F}_i(E,\alpha)\neq\emptyset, \\ \\
0, & \ \text{otherwise}.
\end{array}
\right.$$
\end{prop}

So, to compute, for instance, $\ell_\Oh\left(\dfrac{\I}{\I(\gamma)}\right)$, one takes a chain 
$$m_E=\alpha^0\leq \alpha^1 \leq \cdots \leq \alpha^m=\gamma,$$
and observe that 
$$ 
\ell_\Oh\left(\dfrac{\I}{\I(\gamma)}\right)=\ell_\Oh\left(\dfrac{\I(\alpha^0)}{\I(\gamma)}\right)=\sum_{j=1}^m \ell_\Oh\left(\dfrac{\I(\alpha^{j-1})}{\I(\alpha^j)}\right).
$$

This shows that the colength of $\I(\gamma)$ in $\I$ is independent from the chain one chooses.
So, if we take a saturated chain, that is, a chain such that $\alpha^j-\alpha^{j-1}\in\{e_i; \ i=1,2,\ldots,r\}$ for all $j=1,\ldots,m$, then the computation may be done using Proposition \ref{prop1}.

In what follows We will denote $\ell_\Oh$  simply by $\ell$.




\subsection{Case r=2}

This simplest case was studied by D'anna, Barucci and Fr\"oberg in \cite{BDF1} and we reproduce it here because it gives a clue on how to proceed in general. 

Let $\alpha^0=m_E$ and consider the saturated chain in $\mathbb{Z}^2$ 
$$\alpha^0\leq \cdots\leq\alpha^m=\gamma=(\gamma_1,\gamma_2)$$
such that
$$
\begin{array}{l}
\alpha^0=(\alpha^0_1,\alpha^0_2),\, \alpha^1=(\alpha^0_1+1,\alpha^0_2),\ldots,\alpha^s=(\gamma_1,\alpha^0_2),\\ \\
\alpha^{s+1}=(\gamma_1,\alpha^0_2+1),\, \alpha^{s+2}=(\gamma_1,\alpha^0_2+2) , \ldots,\alpha^m=(\gamma_1,\gamma_2),
\end{array}$$
and consider the following sets
$$L_1=\{\alpha^0,\alpha^1,\ldots, \alpha^s\} \ \ \text{and} \ \ L_2=\{\alpha^s,\alpha^{s+1},\ldots, \alpha^m\}.
$$

By Proposition \ref{prop1}, we have
\begin{eqnarray}
\ell\left(\dfrac{\I}{\I(\gamma)}\right)&=&\#L_1-\#\{\alpha\in L_1; \ \overline{F}_1(E,\alpha)=\emptyset\}+\nonumber \\
&&\#L_2-\#\{\alpha\in L_2; \ \overline{F}_2(E,\alpha)=\emptyset\}.\nonumber
\end{eqnarray}

Now, because of our choice of $L_1$, we have that 
$$
\forall \, \alpha\in L_1, \ \overline{F}_1(E,\alpha)=\emptyset \ \Longleftrightarrow \ pr_1(\alpha)\in M(E_1),$$
hence 
$$\#\{\alpha\in L_1; \ \overline{F}_1(E,\alpha)=\emptyset\}=\#M(E_1).$$

Observe that not all $\alpha\in L_2$ with $\overline{F}_2(E,\alpha)=\emptyset$ are such that $pr_2(\alpha)\in M(E_2)$, hence
$$\#\{\alpha\in L_2; \ \overline{F}_2(E,\alpha)=\emptyset\}=\#M(E_2)-\xi,$$
where $\xi$ is the number of $\alpha$  in $L_2$  with $pr_2(\alpha)\in E_2$ and $\overline{F}_2(E,\alpha)=\emptyset$. But such $\alpha$ are in one-to-one correspondence with the maximals of $E$, hence $\xi=\#M(E)$.\medskip

Putting all this together, we get
\begin{prop} We have that
\begin{equation}\label{I}
\ell\left(\dfrac{\I}{\I(\gamma)}\right)=(\gamma_1-\alpha^0_1)-\#M(E_1) +(\gamma_2-\alpha^0_2)-\#M(E_2)-\#M(E).
\end{equation}
\end{prop}

\subsection{Case $r\geq 3$}

Let us assume that $\I$ is a fractional ideal of $\Oh$, where $\Oh$ has $r$ minimal primes.

\begin{figure}[H]
 \centering
\includegraphics[scale=0.4]{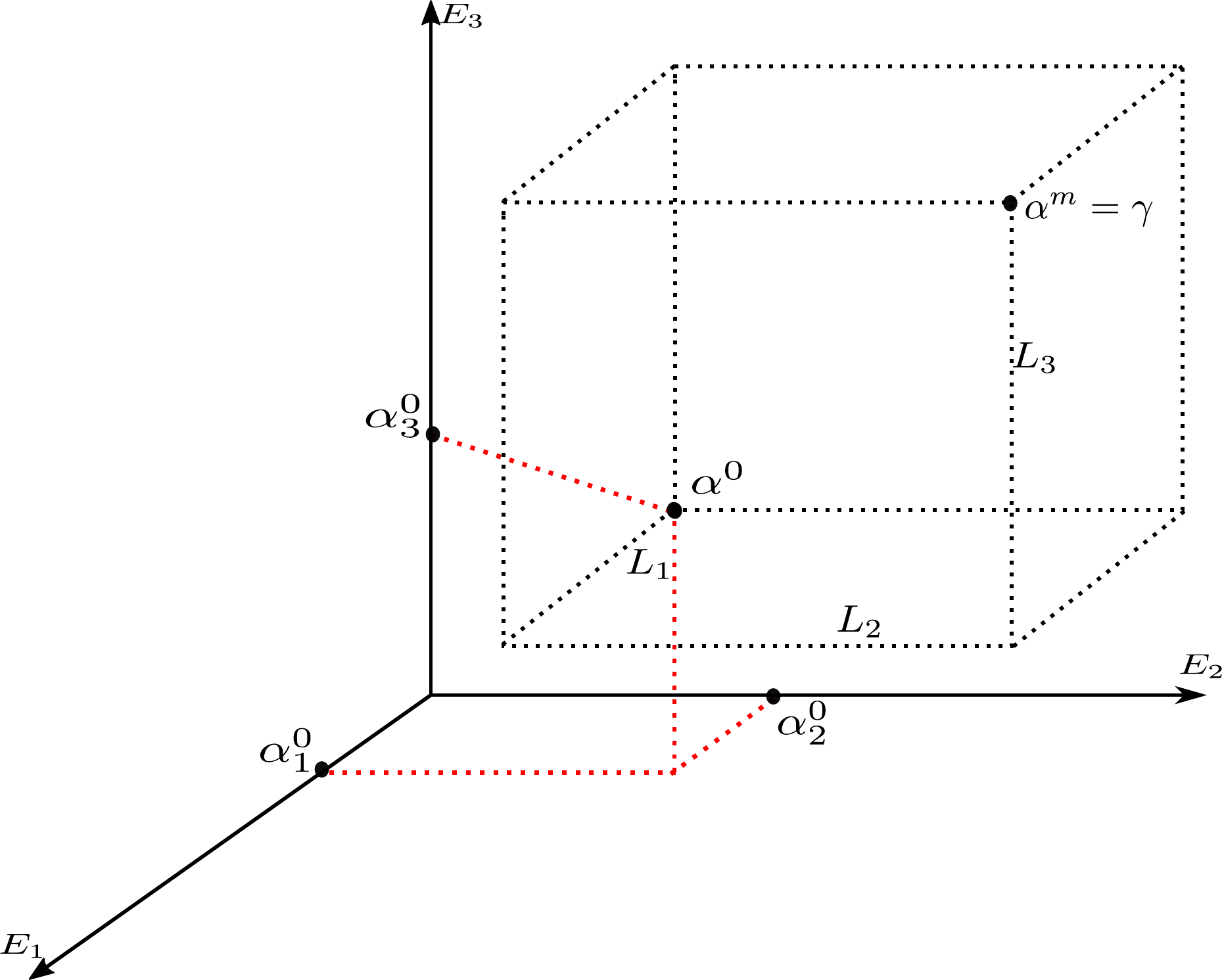}
\caption{The saturated chain for $r=3$}
\label{fig2}
\end{figure}

Let  
 $$m_E=\alpha^0\leq \alpha^1\leq \cdots\leq \alpha^m=\gamma,$$
be the saturated chain in $\mathbb{Z}^r$, given by the union of the following paths (see  \ref{fig2}, for $r=3$):
$$
\begin{array}{l}
L_1\colon \alpha^0, \alpha^1= \alpha^0+e_1,\ldots,\alpha^{s_1}=\alpha^0+(\gamma_1-\alpha^0_1)e_1=(\gamma_1,\ldots,\alpha^0_r),\\ \\
\ldots \\ \\
L_r\colon \alpha^{s_{r-1}}=(\gamma_1,\ldots,\gamma_{r-1},\alpha^0_r), \alpha^{s_{r-1}+1}= \alpha^{s_{r-1}}+e_r, \ldots, \alpha^{m}=
\gamma.
\end{array}
$$

Let us define $I'=\{1,\ldots,r-1\}$. We will need the following result:

\begin{lema}\label{claim} For any $\alpha\in L_1\cup\ldots \cup L_{r-1}$, and for $i\in I'=\{1,\ldots,r-1\}$, one has
$$\overline{F}_i(E,\alpha)\neq\emptyset \ \Longleftrightarrow \ \overline{F}_i(E_{I'},pr_{I'}(\alpha))\neq\emptyset.$$ 
\end{lema}
\pf ($\Rightarrow$) \ This is obvious.

\noindent ($\Leftarrow$)  \ Suppose that $$(\theta_1,\ldots,\theta_{r-1})\in \overline{F}_i(E_{I'}, pr_{I'}(\alpha))\neq\emptyset.$$  Since by Lemma \ref{proj} one has that $pr_{I'}(E)=E_{I'}$, then there exists $\theta=(\theta_1,\ldots,\theta_{r-1},\theta_r)\in E$. Since $\alpha\in L_i$ for some $i=1,\ldots,r-1$, it follows that $\alpha_r=\alpha^0_r$. Then one cannot have $\theta_r<\alpha_r=\alpha^0_r$, because otherwise $$(\alpha^0_1,\ldots,\alpha^0_{r-1},\theta_r)=\min(\alpha^0,\theta)\in E,$$
which is contradiction, since  $\alpha^0$ is the minimum of $E$. Hence $\theta_r\geq\alpha_r$, so $\theta \in \overline{F}_i(E,\alpha)$, and the result follows. 
\cqd

Lemma \ref{claim} allows us to  write:
\begin{equation}\label{intermediateformula}
\ell\left(\dfrac{\I}{\I(\gamma)}\right)=\ell\left(\dfrac{\pi_{I'}(\I)}{\pi_{I'}(\I)(pr_{I'}(\gamma))}\right)+(\gamma_r-\alpha^0_r)-
                                       \#\{\alpha\in L_r; \ \overline{F}_r(E,\alpha)=\emptyset\}.
\end{equation}

Hence to get an inductive formula for $\ell\left(\frac{\I}{\I(\gamma)}\right)$, we only have to compute $$\#\{\alpha\in L_r; \ \overline{F}_r(E,\alpha)=\emptyset\},$$ and for this we will need the following lemma.

\begin{lema}\label{fibra}
Let $\alpha\in\mathbb{Z}^r$, then $\overline{F}_j(E,\alpha)=\emptyset$ if only if there exist some $J\subseteq I$ with $j\in J$ and a relative maximal $\beta$ of $E_J$ such that $pr_{j}(\beta)=\alpha_j$ and $pr_i(\beta)<\alpha_i$, for all $i\in J$, $i\neq j$.
\end{lema}
\pf ($\Leftarrow$) \ (We prove more, since it is enough to assume $\beta$ is any maximal) Let us assume that there exist $ J\subseteq I$, with $j\in J$ and $\beta\in M(E_J)$, such that $\beta_j=pr_{j}(\beta)=\alpha_j$ and $\beta_i=pr_i(\beta)<\alpha_i$, for all $i\in J$, $i\neq j$. 

Suppose by reductio ad absurdum that  $\overline{F}_j(E,\alpha)\neq\emptyset$. Let $\theta\in\overline{F}_j(E,\alpha)$ that is $\theta_j=\alpha_j$ and $\theta_i\geq\alpha_i, \forall i\in J\setminus \{j\}$. Now since 
$$pr_j(pr_J(\theta))=\theta_j=\alpha_j=\beta_j \ \text{and} \ pr_i(pr_J(\theta))=\theta_j\geq \alpha_i>\beta_i, \forall i\in J, i\neq j,$$
then $pr_J(\theta)\in F_j(E_J,\beta)$, which contradicts the assumption that  $\beta\in M(E_J)$.  \smallskip
	
\noindent ($\Rightarrow$)  \ Since $\overline{F}_j(E,\alpha)=\emptyset$ implies $F_j(E,\alpha)=\emptyset$, the proof follows the same lines as the proof of \cite[Theorem 1.5]{D88}.\cqd

Going back to our main calculation, by Lemma \ref{fibra}, if $\alpha\in L_r$ is such that $\overline{F}_r(E,\alpha)=\emptyset$, then there exist a subset $J$ of  $I=\{1,\ldots,r\}$, with $r\in J$, and $\beta\in RM(E_J)$, with $pr_r(\beta)=\alpha_r$ and $pr_i(\beta)<\alpha_i$ for $i\in J, i\neq r$. 

Notice that for $\alpha \in L_r$ one has $\alpha_i=\gamma_i$ for $i\neq r$, so the condition $pr_i(\beta)<\alpha_i$ for $i\in J, i\neq r$ is satisfied, since $\beta\in M(E_J)$. So, we have a bijection between the set $\{\alpha\in L_r; \ \overline{F}_r(E,\alpha)=\emptyset\}$ and the set $\bigcup_{r\in J\subseteq I}pr_r(RM(E_J))$, where $pr_r(\beta)$ is the  last coordinate of $\beta$, that is, $\beta_r$.

This gives us the formula

$$\#\{\alpha\in L_r; \ \overline{F}_r(E,\alpha)=\emptyset\}=\#\bigg(\bigcup_{r\in J\subseteq I}pr_r(RM(E_J))\bigg).
$$

If $J=\{r\}$, we have $pr_r(M(E_r))= M(E_r)$, and for all $J'\subseteq I$ with $r\in J'$ and $J'\neq\{r\}$, also $M(E_r)\cap pr_r(RM(E_{J'}))=\emptyset$,  hence
\begin{equation} \label{formula4}\#\{\alpha\in L_r;\overline{F}_r(E,\alpha)=\emptyset\}=\#M(E_r)+\#\bigg(\bigcup_{\stackrel{ r\in J\subseteq I}{J\neq\{r\}}}pr_r(RM(E_J))\bigg). 
\end{equation}

Now, putting together Equations (\ref{intermediateformula}) and (\ref{formula4}), we get the following recursive formula: 

\begin{teo} \label{teogeral1} For a fractional ideal $\I$ of a ring $\Oh$ with $r$ minimal primes, with values set $E$, one has
\begin{equation} \label{formulafinal1}\begin{array}{rcl}
\ell\left(\dfrac{\I }{\I(\gamma)}\right)=\ell\left(\dfrac{\pi_{I'}(\I)}{\pi_{I'}(\I)(pr_{I'}(\gamma))}\right)&+&(\gamma_r-\alpha^0_r)-\#M(E_r)-\\
& & \#\bigg(\bigcup_{\stackrel{ r\in J\subseteq I}{J\neq\{r\}}}pr_r(RM(E_J))\bigg).
\end{array}
\end{equation}
\end{teo}





\section{Another characterization for Gorensteiness}

The central results in this section  will be other two  characterizations of the Gorensteiness of the ring of an algebroid curve, namely Theorems \ref{dimGo} and \ref{p+q} which are respectively generalizations of \cite[Theorems 3.6 and 5.3]{CDK}. As a consequence, we get a corollary that establishes a correspondence between relative maximals of the set of values of a fractional ideal of a Gorenstein algebroid curve and the absolute maximals of the set of values of its dual ideal, generalizing the symmetry relation among relative maximals and absolute maximals in the semigroup of values of a Gorenstein ring, described by F. Delgado in \cite{D88} for plane algebroid curves and by Campillo, Delgado and Kyiek in \cite{CDK} for algebroid Gorenstein curves. This last referred correspondence was also recently obtained by D. Pol in \cite{Po2} by using other methods. 

\begin{prop}\label{l_E}
	Let $\I$ be a fractional ideal of $\Oh$ and $\alpha\in\mathbb{Z}^r$. Then
	$$\ell\left(\frac{\I(\alpha)}{\I(\alpha+e_i)}\right)+\ell\left(\frac{\I^*(c(S)-\alpha-e_i)}{\I^*(c(S)-\alpha)}\right)\leq1 \ \text{for every} \ i\in I.$$
	\end{prop}
\pf
By Proposition \ref{prop1} we have $\ell\left(\frac{\I(\alpha)}{\I(\alpha+e_i)}\right)\leq1$ and $\ell\left(\frac{\I^*(c(S)-\alpha-e_i)}{\I^*(c(S)-\alpha)}\right)\leq1$. Thus, it is enough to show that if $\ell\left(\frac{\I(\alpha)}{\I(\alpha+e_i)}\right)=1$ then $\ell\left(\frac{\I^*(c(S)-\alpha-e_i)}{\I^*(c(S)-\alpha)}\right)=0$ and vice versa.

Proposition \ref{prop1} also says that $\ell\left(\frac{\I(\alpha)}{\I(\alpha+e_i)}\right)=1$ if and only if $\overline{F}_i(E,\alpha)\neq\emptyset$ and that $\ell\left(\frac{\I^*(c(S)-\alpha-e_i)}{\I^*(c(S)-\alpha)}\right)=0$ if and only if $\overline{F}_i(E^*,c(S)-\alpha-e_i)=\emptyset.$ 

Let us assume by reductio ad absurdum that $$\overline{F}_i(E,\alpha)\neq\emptyset \qquad \text{and} \qquad \overline{F}_i(E^*,c(S)-\alpha-e_i)\neq \emptyset.$$ Take $\theta$ in the first of these sets and $\theta'$ in the second, then $\theta+\theta'\in S$; even more, we have that 
 $\theta+\theta'\in F_i(S,\f(S))$ because $\theta_i+\theta'_i=\f_i(S)$ and $\theta_j+\theta'_j>\f_j(S)$ for all $j\neq i$, hence $\theta+\theta'\in F_i(S,\f(S))$, which is a contradiction, since $F_i(S,\f(S))=\emptyset$.

The proof of the other case is similar.
\cqd

The following theorem generalizes \cite[Theorem 3.6]{CDK}.
\begin{teo}\label{dimGo}
	Let $\Oh$ be a finitely generated local reduced and complete one-dimensional algebra. Then the following statements are equivalent:
	\begin{enumerate}[\rm i)]
	\item For all fractional ideal $\I$ of $\Oh$, one has 
		
	$\ell\left(\frac{\I(\alpha)}{\I(\alpha+e_i)}\right)+\ell\left(\frac{\I^*(c(S)-\alpha-e_i)}{\I^*(c(S)-\alpha)}\right)=1, \ \text{for every} \ \alpha\in\mathbb{Z}^r \text{and} \ i\in I.$
	\item $\Oh$ is a Gorenstein ring.
	\end{enumerate}	
\end{teo}
\pf
\  ${\rm i) \Rightarrow ii)}$\  \
Assuming (i) and taking $\I=\Oh$, we will show that $S$ is symmetric, proving that $\Oh$ is Gorenstein.

Recall that $S$ symmetric means that 
\[
\alpha\in S \ \Longleftrightarrow \ F(S,\f(S)-\alpha)=\emptyset.
\]

One direction of the above equivalence is known to be true from Corollary \ref{halfsymm}, namely
\[
\alpha\in S \ \Rightarrow \ F(S,\f(S)-\alpha)=\emptyset.
\]

To prove the other direction, we show firstly that 
\begin{equation} \label{eq}
F(S,\f(S)-\alpha)=\emptyset\Rightarrow\overline{F}_i(S,c(S)-\alpha-e_i)=\emptyset, \ \text{for all} \ i\in I.
\end{equation}

Indeed, if $F(S,\f(S)-\alpha)=\emptyset$ and $\overline{F}_i(S,c(S)-\alpha-e_i)\neq\emptyset$, for some $i\in I$, take 
$\theta \in \overline{F}_i(S,c(S)-\alpha-e_i)$, hence $\theta_i=c_i(S)-\alpha_i-1=\f_i(S)-\alpha$ and $\theta_j\geq c_j(S)-\alpha_j>\f_i(S)-\alpha_j$, for all $j\neq i$. This implies that $\theta\in F_i(S,\f(S)-\alpha)\subset F(S,\f(S)-\alpha)=\emptyset$, a contradiction. So, we proved (\ref{eq}).

From (\ref{eq}) and Proposition  \ref{prop1}, we get that $\ell\left(\frac{\Oh(c(S)-\alpha-e_i)}{\Oh(c(S)-\alpha)}\right)=0$, for all $i\in I$. This in view of (i) implies that $\ell\left(\frac{\Oh(\alpha)}{\Oh(\alpha+e_i)}\right)=1$, for all $i\in I$, which in turn, in view of Proposition \ref{prop1} implies that $\overline{F}_i(S,\alpha)\neq\emptyset$ for all $i\in I$, from which we conclude that $\alpha\in S$, as we wanted to show.

\noindent ${\rm ii) \Rightarrow i)}$ \ \ We suppose now that $\Oh$ is Gorenstein and $\I$ is a fractional ideal of $\Oh$. 

From Theorem \ref{poltheo} we know that 
\[
\alpha\notin E \ \Leftrightarrow \ F(E^*,\f(S)-\alpha)\neq \emptyset.
\]

To conclude the proof of this part of the theorem, it is enough to show that 
\[
\forall \ i\in I, \ \ \overline{F}_i(E,\alpha)=\emptyset \ \Longrightarrow \ \overline{F}_i(E^*,\f(S)-\alpha)\neq \emptyset .
\]

Indeed, because the converse of the above implication is always true (cf. Proposition \ref{l_E}), item (i) will follow by Proposition \ref{prop1}.

Suppose now that $\overline{F}_i(E,\alpha)=\emptyset$. This implies that  $\alpha\notin E$, so, in view of Theorem \ref{poltheo} we get that $F(E^*,\f(S)-\alpha)\neq \emptyset$. Hence there exists an element $\theta \in F_i(E^*,\f(S)-\alpha)$, for some $i\in I$. So, $\theta_i=\f_i(S)-\alpha_i=c(S)-1-\alpha_i$ and $\theta_j>\f_j(S)-\alpha_j=c_j(S)-1-\alpha_j$, for all $j\neq i$. This last inequlity implies that $\theta_j\geq c_j(S)-\alpha_j$, for all $j\neq i$. Thus $\theta\in \overline{F}_i(E^*,c(S)-\alpha-e_i)$, hence this set is nonempty.
\cqd

In \cite[Section 5]{CDK}, the authors define certain integers related to the semigroup of values of an algebroid curve to give another characterization for the property of $\Oh$ being Gorenstein. Now, let us define these numbers for fractional ideals.\medskip

We denote $\ell\left(\frac{\I(\alpha)}{\I(\alpha+e_i)}\right)$ by $\ell_E(\alpha,e_i)$ and put $e_J=\sum_{i\in J}e_i$ for $\emptyset \neq J\subset I$, $e_\emptyset=0$ and $e=e_I=(1,1,\ldots,1)$. We will also use the notation $J^c=I\setminus J$.

Let us define the following sentences:
$$
p(\alpha,j) \ : \ \forall J\subset I, \ \text{with} \ \#J=j, \quad \ell_E(\alpha+e_{J^c},e_i)=0, \ \forall i\in J, \ \ 1\leq j\leq r.$$

Notice that $p(\alpha,n)\Rightarrow p(\alpha,n-1)$. \smallskip

Indeed, let $K\subset I$ with $\#K=n-1$, now $K$ is included in some subset $J$ of $I$ with $\#J=n$. From the assertion $p(\alpha,n)$ we have that $\ell_E(\alpha+e_{J^c},e_i)=0$ for all $i\in J$. Suppose by reductio ad absrurdum that $p(\alpha,n-1)$ is false, hence 
$\ell_E(\alpha+e_{K^c},e_i)=1$ for some $i\in K$, that is, $\overline{F}_i(E,\alpha+e_{K^c})\neq\emptyset$. Take $\theta \in \overline{F}_i(E,\alpha+e_{K^c})$, then $\theta\in\overline{F}_i(E,\alpha+e_{J^c})$, because $\theta_i=\alpha_i$ and $\theta_k\geq\alpha_k$ for $k\in K$, $\theta_j\geq\alpha_j+1$ for $j\in K^c$. Hence $\overline{F}_i(E,\alpha+e_{J^c})\neq \emptyset$, which contradicts $\ell_E(\alpha+e_{J^c},e_i)=0$ for all $i\in J$.
\medskip

Thus, we can define
$$p(E,\alpha)=\max\{n; \ell_E(\alpha+e_{J^c},e_i)=0, \forall J\subset I, \#J \leq n, \forall i\in J\}.$$

By the definition of $p(E,\alpha)$ it follows that 
\begin{enumerate}[i)]
	\item $0\leq p(E,\alpha)\leq r$, for every $\alpha$;
	\item $p(E,\alpha)=\max\{n; F_A(E,\alpha)=\emptyset,\forall A\subset I,\#A\leq n\}$;
		\item $p(E,\alpha)\geq 1\Longleftrightarrow F(E,\alpha)=\emptyset$.
\end{enumerate}

We also define the sentences:
$$\begin{array}{l}
q(\alpha,j)\ : \ \forall J\subset I, \ \text{with} \ \#J=j, \quad \ell_E(\alpha+e_{J^c},e_i)=1, \forall i\in J, \ 1\leq j\leq r;\\ \\
q(\alpha, r+1) \ : \ F(E,\alpha)=\emptyset.
\end{array}
$$

Notice that $q(\alpha,n-1)\Rightarrow q(\alpha,n)$. \smallskip

Indeed, assume that $q(\alpha,n-1)$ is true. Let $K\subset I$, $\#K=n$ and let $J$ be a subset of $K$ with $\#J=n-1$,  we have $\ell_E(\alpha+e_{J^c},e_i)=1$ for all $i\in J$, i.e., $\overline{F}_i(E,\alpha+e_{J^c})\neq\emptyset$. Take $\theta \in \overline{F}_i(E,\alpha+e_{J^c})$ and notice that $\theta\in\overline{F}_i(E,\alpha+e_{K^c})$, because $\theta_i=\alpha_i$ and $\theta_j\geq\alpha_j$ for $j\in K$, $\theta_j\geq\alpha_j+1$ for $j\in K^c$, hence $\ell_E(\alpha+e_{K^c},e_i)=1$, for all $i\in K$.\medskip

Thus, we can also define
$$q(E,\alpha)=\min\{n; \ell_E(\alpha+e_{J^c},e_i)=1, \forall J\subset I, \#J\geq  n, \forall i\in J\}.$$

By definition $q(E,\alpha)$ has that following properties;
\begin{enumerate}[i)]
	\item $1\leq q(E,\alpha)\leq r+1$, for every $\alpha$;
	\item $q(E,\alpha)=\min\{n; F_A(E,\alpha)\neq\emptyset,\forall A\subset I,\#A\geq n\}$
	\item $q(E,\alpha)\leq r\Longleftrightarrow \overline{F}_i(E,\alpha)\neq\emptyset$, for $i=1,\ldots,r$;
	\item $p(E,\alpha)<q(E,\alpha)$ for every $\alpha$.
\end{enumerate}

The following theorem generalizes \cite[Theorem 5.3]{CDK}.
\begin{teo}\label{p+q}
	With notations as above, for any fractional ideal $\I$ of $\Oh$ with value set $E$ and for any $\alpha\in\mathbb{Z}^r$, we have
\begin{equation}\label{p,q}
p(E,\alpha)+q(E^*,\f(S)-\alpha)\geq r+1.
\end{equation}
Moreover, $\Oh$ is a Gorenstein ring if an only if equality holds in (\ref{p,q}) for every fractional ideal $\I$ of $\Oh$ and every $\alpha\in\mathbb{Z}^r$.
\end{teo}
\pf
Take any fractional ideal $\I$ and $\alpha\in\mathbb{Z}^r$, and define $n=r+1-q(E^*,f(S)-\alpha)$. Take $J\subset I$ with $\#J\leq n$ and any $i\in J$. Put $K=J^c\cup \{i\}$. By Proposition \ref{l_E}, one has
\begin{equation}\label{eq'}
\ell_E(\alpha+e_{J^c},e_i)+\ell_{E^*}(c(S)-e-\alpha+e_{K^c},e_i)\leq1, \ \text{for all} \ i\in J.
\end{equation}

Since  $\#K\geq q(E^*,f(S)-\alpha)$ and $i\notin K^c$, by definition of $q$ one has $\ell_{E^*}(c(S)-e-\alpha-e_{K^c},e_i)=1$ and therefore $\ell_E(\alpha+e_{J^c},e_i)=0$. This shows that $n\leq p(E,\alpha)$, and therefore (\ref{p,q}) is proven. 

Now, if $\Oh$ is Gorenstein, to prove that (\ref{p,q}) holds, we must show that $p(E,\alpha)=n$. Suppose by reductio ad absurdum that 
$p(E,\alpha)\geq n+1$. Then, we have that $$\ell_E(\alpha+e_{J^c},e_i)=0, \ \forall i\in J, \ \forall J\subset I \ \text{with}\ \#J=n+1.$$
Consider a set $K\subset I$ with $\#K=r-n$ and define the set $J=K^c\cup \{i\}$, with $i\in K$. Notice that $J$ has $n+1$ elements, then
 $$\ell_E(\alpha+e_{J^c},e_i)=0, \ \forall i\in J.$$
Since, that $\Oh$ is Gorenstein, we have
$$\ell_{E*}(\f(S)-\alpha+e_{K^c},e_i)=0,\ \forall i\in K,$$
and since, $i$ is any element of $K$, we have that
$$\ell_{E*}(\f(S)-\alpha+e_{K^c},e_i)=0,\ \forall i\in K \ \text{and} \ \forall K\subset I, \ \text{with}\ \#K=r-n.$$
As $n=r+1-q(\f(S)-\alpha)$, we have that $\#K<q(\f(S)-\alpha)$, this contradicts the definition of $q(\f(S)-\alpha)$. Therefore, it follows that $p(E,\alpha)=n$.

Now, assume that we have equality in (\ref{p,q}), for every fractional ideal $\I$ of $\Oh$ and every $\alpha\in\mathbb{Z}^r$. Let $\alpha\in\mathbb{Z}^r$ and $i\in I$. If $\overline{F}_i(E,\alpha)=\emptyset$, then there exists $\beta$ with $\beta_i=\alpha_i$ and $\beta_j<\alpha_j$ for every $j\neq i$ such that $F(E,\beta)=\emptyset$. From the condition $F(E,\beta)=\emptyset$ we get that $p(E,\beta)\geq1$, hence from Equality (\ref{p,q}) we get that  $q(E^*,\f(S)-\beta)\leq r$. Therefore, $\overline{F}_j(E^*,\f(S)-\beta)\neq\emptyset$ for every $j\in I$. In particular, $\overline{F}_i(E^*,\f(S)-\beta)\neq\emptyset$. On the other hand, since $\f_i(S)-\beta_i=\f_i(S)-\alpha_i$ and $\f_j(S)-\beta_j> \f_j(S)-\alpha_j$, it follows that $\f(S)-\beta\geq\f(S)-\alpha$, therefore $\overline{F}_i(E^*,\f(S)-\alpha)\neq\emptyset$. So, we proved that
$$\ell_E(\alpha,e_i)+\ell_{E^*}(c(S)-\alpha-e_i,e_i)=1, \ \text{for every} \ \alpha\in\mathbb{Z}^r \text{and} \ \text{all} \ i\in I.$$

So, by Theorem \ref{dimGo} it follows that $\Oh$ is Gorenstein.
\cqd

Notice that in Theorem \ref{p+q} we have the same result with the inequality
\[
p(E^*,\beta)+q(E,\f(S)-\beta)\geq r+1.
\]

\begin{rem}
In terms of $p(E,\alpha)$ and $q(E,\alpha)$, we may characterize the maximals of $E$ as follows:\medskip

\noindent {\rm i)} $\alpha$ is maximal if and only if $p(E,\alpha)\geq1$ and $q(E,\alpha)\leq r$; \medskip

\noindent {\rm ii)} $\alpha$ is a relative maximal if and only if $p(E,\alpha)=1$ and $q(E,\alpha)=2$; \ and \medskip

\noindent {\rm iii)} $\alpha$ is an absolute maximal if and only if $p(E,\alpha)=r-1$ and $q(E,\alpha)=r$.
\end{rem}

The following result will generalize \cite[Theorem 2.10]{D87}, and is a consequence of Theorem \ref{p+q} and the previous remark.

\begin{cor}[\textbf{Symmetry of maximals}]\label{ed3} Suppose that $\mathcal{O}$ is Gorenstein with value set $S$ and $\I$ is a fractional ideal of $\Oh$ with value set $E$. Then for  $\beta\in E^*$, one has that $\beta$ is a maximal of $E^*$ if and only if $\f(S)-\beta\in E$. Moreover, if $\beta\in E^*$ and $\alpha\in E$ are such that $\alpha+\beta=\f(S)$, then $\alpha$ is an absolute maximal of $E$ if and only if $\beta$ is a relative maximal of $E^*$.
\end{cor}
\pf
Suppose that $ \beta $ is a maximal of $ E^*$, then

$$F(E^*,\f(S)-(\f(S)-\beta))=F(E^*,\beta)=\emptyset,$$
hence from (\ref{pol}) it follows that $\f(S)-\beta\in E $.

Conversely, suppose that $\f(S)-\beta\in E$, then from (\ref{pol}) we have that $F(E^*,\beta)=\emptyset$ and, as by hypothesis, $\beta\in E^*$, it follows that $\beta$ is a maximal of $E^*$. This proves the first part of the theorem.\smallskip

We now prove the second part of the theorem. Suppose that $\beta$ is a relative maximal of $E^*$ and that $\alpha\in E$ is such that $\alpha+\beta=\f(S)$. From the first part of the theorem that we have just proved, it follows that $\alpha$ is a maximal of $E$.

By reductio ad absurdum, suppose that $\alpha$ is not an absolute maximal of $E$, then there exists $J\subset I$ with $2\leq \#J\leq r-1$ such that $F_J(E,\alpha)\neq\emptyset$. On the other hand, since $\beta$ is a relative maximal of $E^*$, then $F_K(E^*,\beta)\neq\emptyset$, for all $K\subset I$ with $\#K\geq 2$. In particular, take  $K=(I\setminus J)\cup\{i\}$, where $i\in J$, and choose $\theta\in F_K(E^*,\beta)$ and $\theta'\in F_J(E,\alpha)$. Notice that because of Equation (\ref{$E+E^*$}) one has $\theta+\theta'\in S$. On the other hand, our choice of $K$ implies that $\theta+\theta'\in F(\Z^r,\f(S))$, hence $\theta+\theta'\in F(S,\f(S))$, which is a contradiction, since from Corollary \ref{halfsymm} we know that $F(S,\f(S))=\emptyset$. Thus $\alpha$ is an absolute maximal of $E$.

It remains to prove that if $\alpha+\beta=\f(S)$ and $\alpha$ is an absolute maximal of $E$, then $\beta$ is a relative maximal of $E^*$. In fact, suppose that $\alpha$ is an absolute maximal of $E$, which means that $p(E,\alpha)=r-1$ and $q(E,\alpha)=r$, then by Theorem \ref{p+q} we have that $q(E^*,\beta)=2$ and $p(E^*,\beta)=1$, i.e., $\beta$ is a relative maximal of $E^*$.
\cqd

\end{document}